\documentclass{amsproc}

\usepackage{latexsym}
\usepackage{amsmath}
\usepackage{amsfonts}
\usepackage{amssymb}
\usepackage{amscd}

\numberwithin{equation}{section}

\newtheorem{proposition}{Proposition}[section]
\newtheorem{lemma}[proposition]{Lemma}
\newtheorem{theorem}[proposition]{Theorem}
\newtheorem{corollary}[proposition]{Corollary}
\newtheorem{definition}[proposition]{Definition}

\newtheorem{remark}[proposition]{Remark}

\def\itheorem#1#2{\newtheorem{#1}[proposition]{#2}}

\itheorem{Remark}{Remark}

\def\Ext{\mathop{\rm Ext}\nolimits}

\def\rg{\mathop{\rm rg}\nolimits}

\def\cf{\mathop{\rm cf}\nolimits}

\def\Q{{\mathbb{Q}}}
\def\Z{{\mathbb{Z}}}

\def\C{{\mathfrak{C}}}

\def\tc{{\mathcal{T}\mathcal{C}}}

\def\P{{\mathcal{P}}}

\begin{document}

\title{Some ugly $\aleph_1$-free abelian groups}

\author{Saharon Shelah}

\address{The Hebrew University, Givat Ram, Jerusalem 91904,
Israel and \\ Rutgers University, Newbrunswick, NJ U.S.A.}

\email{Shelah@math.huji.ae.il}

\thanks{Publication 773 in the first author's list of publication. The first author was supported
by project No. G-0545-173,06/97 of the {\em German-Israeli
Foundation for Scientific Research \& Development.}}

\author{Lutz Str\"ungmann}

\address{Fachbereich 6 -- Mathematik,
University of Essen, 45117 Essen, Germany}

\curraddr{The Hebrew University, Givat Ram, Jerusalem 91904,
Israel}

\email{lutz@math.huji.ac.il}

\thanks{The second author was supported by a MINERVA fellowship.}


\begin{abstract}
Given an $\aleph_1$-free abelian group $G$ we characterize the
class $\C_G$ of all torsion abelian groups $T$ satisfying
$\Ext(G,T)=0$ assuming the special continuum hypothesis $CH$.
Moreover, in G\"odel's constructable universe we prove that this
characterizes $\C_G$ for arbitrary torsion-free abelian $G$. It
follows that there exist some ugly $\aleph_1$-free abelian groups.
\end{abstract}

\maketitle

\section{Introduction}
In 1969 Griffith \cite{G} solved Baer's splitting problem on mixed
abelian groups when he proved that an abelian group $G$ is free if
and only if $\Ext(G,T)=0$ for all torsion abelian groups $T$. It
is easy to see that an abelian group $G$ which satisfies
$\Ext(G,T)=0$ for all torsion abelian groups $T$ must be
torsion-free and homogeneous of type $\Z$. Thus it was natural to
ask whether or not one could extend Griffith's result to
homogeneous torsion-free groups which are not necessarily of
idempotent type. That this is not the case was shown in
\cite{St1} by the second author. This was a consequence of
techniques and results obtained in \cite{SW}. Inspired by Baer's
question \cite{B} to characterize all pairs of torsion-free
abelian $G$ and torsion abelian $T$ such that $\Ext(G,T)=0$,
Wallutis and the second author considered in \cite{SW} the
torsion groups of the cotorsion class singly cogenerated by a
torsion-free group $G$. Cotorsion theories were introduced by
Salce in \cite{S} but it was the first time in \cite{SW} that
only the torsion groups of the cotorsion theory were considered.
Recall, that for a torsion-free abelian group $G$ the class of all
torsion abelian groups $T$ satisfying $\Ext(G,T)=0$ is denoted by
$\tc(G)$ (see \cite{SW}). This class is obviously closed under
taking epimorphic images and contains all torsion cotorsion
groups, i.e. all bounded groups. In \cite{SW} satisfactory
characterizations of $\tc(G)$ were obtained for countable
torsion-free abelian groups and for completely decomposable
groups. In fact, it was proved in \cite{SW} that for every
countable torsion-free abelian group $G$ there exists a
completely decomposable group $C$ such that $\tc(G)=\tc(C)$. It
was later shown in \cite{St1} by the second author that for every
finite rank torsion-free abelian group $G$ there even exists a
rational group $R\subseteq \Q$ such that $\tc(G)=\tc(R)$. Thus,
knowing the class $\tc(C)$ for completely decomposable groups
$C$, it was reasonable to search for groups $G$ of uncountable
cardinality such that $\tc(G)$ equals $\tc(C)$ for some completely
decomposable group $C$. Although a criterion was found in \cite[Theorem 3.6]{SW} for
characterizing those classes of torsion abelian groups which may
appear as $\tc(C)$ for completely decomposable group $C$, it
remained open if for instance in G\"odel's universe every
torsion-free abelian group is of this kind. It shall be shown in
this paper that this is not the case but it holds if we replace
completely decomposable by $\aleph_1$-free of cardinality
$\aleph_1$. \\
Assuming $CH$ we give in section 2 a construction of
$\aleph_1$-free abelian groups $G$ of size $\aleph_1$ having a
strange class $\tc(G)$. It shall be proved that for every ideal
$I$ in the set of primes (more general in the set of all powers
of primes) containing all finite subsets of the set of primes,
there exists an $\aleph_1$-free abelian group $G$ of size
$\aleph_1$ such that $\bigoplus\limits_{p \in P}\Z(p) \in \tc(G)$
if and only if $P \in I$. It follows in section 3 that in
G\"odel's constructable universe ($V=L$) every torsion-free
abelian group $G$ satisfies $\tc(G)=\tc(H)$ for some
$\aleph_1$-free group $H$ of size $\aleph_1$. Thus we obtain a
characterization of the class $\tc(G)$ for all torsion-free
abelian groups $G$ in G\"odel's universe and prove that the
structure of the group $G$ is not very much effected by the class
$\tc(G)$. This solves Baer's problem in $V=L$ and contrasts a
result from \cite{GSW} in which it was shown that the cotorsion
theory singly cogenerated by $G$
determines the group $G$ in its structure.\\

All groups under consideration are abelian. The notations are
standard and for unexplained notions in abelian group theory and
set theory we refer to \cite{Fu} and \cite{EM}.

\section{The construction}

In this section we construct some $\aleph_1$-free abelian groups
having special properties. Let us first recall a definition from
\cite{SW}. For a torsion-free group $G$ we denote by $\tc(G)$ the
class of all torsion groups $T$ satisfying $\Ext(G,T)=0$.
Obviously, the class $\tc(G)$ is closed under taking epimorphic
images and contains all torsion cotorsion groups, i.e. all
bounded groups. Recall that a torsion-free group $G$ is called
{\it $\aleph_1$-free} if all its countable subgroups are free.
Let $\Pi$ be the set of natural primes. By $\bar{\Pi}$ we denote
the set of all powers of natural primes, i.e. $\bar{\Pi}=\{ p^n :
p \in \Pi, n < \omega\}$. Moreover, for an infinite subset $P
\subseteq \bar{\Pi}$ we define $T_P=\bigoplus\limits_{p \in
P}\Z(p)$, where $\Z(p)$ denotes the cyclic group of order $p$.
The reader should keep in mind that here $p$ is not necessarily a
prime but could be a prime power. We begin with a compactness
result for countable torsion-free groups (see \cite[Lemma
3.1]{St1}).

\begin{lemma}[\cite{St1}]
\label{compact} Let $G$ be a countable torsion-free group and $T$
a torsion group. Then $T \in \tc(G)$ if and only if $T \in \tc(H)$
for all finite rank pure subgroups $H$ of $G$.
\end{lemma}

\begin{proof}
The proof can be found in \cite[Lemma 3.1]{St1} and is based on the fact that for countable $G$ and any pure subgroup $H\subseteq G$ of finite rank, $T \in \tc(G)$ implies $T \in \tc(G/H)$ ($T$ a torsion group).
\end{proof}

Recall, that for a torsion-free group $G$ of size $\kappa$ a {\it
$\kappa$-filtration of $G$} is a continuous ascending chain of
pure subgroups of cardinality less than $\kappa$ such
that its union equals $G$.

\begin{proposition}[$CH$]
\label{characterization} Let $G$ be a torsion-free group of
cardinality $\aleph_1$ and $P$ an infinite subset of $\bar{\Pi}$.
If $\left< G_{\alpha} : \alpha < \omega_1 \right>$ is an
$\omega_1$-filtration of $G$, then $T_P \not\in \tc(G)$ if and
only if one of the following conditions holds:
\begin{enumerate}
\item $T_P \not\in \tc(H)$ for some finite rank pure subgroup $H$
of $G$ or;
\item $\{ \delta < \omega_1 : G/G_{\delta} \textit{ contains a
finite rank pure subgroup $L_{\delta} \subseteq G/G_{\delta}$ such
that }\\ T_P \not\in \tc(L_{\delta}) \}$ is stationary in $\omega_1$.
\end{enumerate}
\end{proposition}

\begin{proof}
Let \begin{equation}
 S=\{ \delta < \omega_1 : G/G_{\delta}
\textit{ contains a finite rank pure subgroup $L_{\delta}
\subseteq G/G_{\delta}$} \end{equation}
\[ \textit{ such that }
T_P \not\in \tc(L_{\delta}) \}.\] By Lemma \ref{compact} $S=\{
\delta < \omega_1 : T_P \not\in \tc(G_{\beta}/G_{\delta}) \textit{
for some } \delta < \beta < \omega_1 \}$. Now it is easy to see
that $S$ stationary implies that the relative $\Gamma$-invariant
$\Gamma_{T_P}(G) \not=0$. Since we are assuming $CH$ the weak
diamond $\Phi_{\aleph_1}$ holds (see \cite{DS}). Thus (i) or (ii)
imply $T_P \not\in \tc(G)$ by \cite[Proposition XII.1.15]{EM}.
Conversely, assume that $T_P \not\in \tc(G)$ but (i) and (ii) do
not hold.  Then, the relative $\Gamma$-invariant
$\Gamma_{T_P}(G)=0$ and hence \cite[Theorem XII.1.14]{EM} shows
that $T_P \in \tc(G)$ - a contradiction.
\end{proof}

\begin{remark}
If we assume $V=L$, then we could extend Proposition \ref{characterization} to larger cardinalities using techniques as for instance developed in \cite[Theorem 3.1]{BFS} and using an appropriate filtration but it is not needed here.
\end{remark}
 
Let $S$ be a stationary subset of $\omega_1$ consisting of limit
ordinals, i.e. for all $\alpha \in S$, $\cf(\alpha)=\omega$.
Recall the following definition.

\begin{definition}
A {\it ladder system} $\bar{\eta}$ on $S$ is a family of functions
$\bar{\eta}=\left< \eta_{\delta} : \delta \in S \right>$ such that
$\eta_{\delta}: \omega \rightarrow \delta$ is strictly increasing
with $\sup(\rg(\eta_{\delta}))=\delta$, where $\rg(\eta_{\delta})$
denotes the range of $\eta_{\delta}$. We call the ladder system
{\it tree-like} if for all $\delta, \nu \in \bar{\eta}$ and every
$\alpha, \beta \in \omega$,
$\eta_{\delta}(\alpha)=\eta_{\nu}(\beta)$ implies $\alpha=\beta$
and $\eta_{\delta}(\rho)=\eta_{\nu}(\rho)$ for all $\rho \leq
\alpha$.
\end{definition}

\begin{proposition}[$CH$]
\label{existence} Let $\left< P_{\alpha} : \alpha < \omega_1
\right>$ be a sequence of infinite subsets of $\bar{\Pi}$. Then
there exists an $\aleph_1$-free torsion-free group $G$ of cardinality
$\aleph_1$ such that for any infinite subset $P$ of $\bar{\Pi}$,
$T_P \not\in \tc(G)$ if and only if $\{ \delta < \omega_1 : |P \cap
P_{\delta}|=\aleph_0 \}$ is stationary.
\end{proposition}

\begin{proof}
Since we are assuming $CH$ the weak diamond $\Phi_{\aleph_1}$
holds. Let $S$ be a stationary subset of $\omega_1$ such that
$\Phi_{\aleph_1}(S)$ holds. Since $lim(\omega_1)$ is a cub in
$\omega_1$ we may assume without loss of generality that
$S=lim(\omega_1)$, i.e. $S$ consists of all limit ordinals of
$\omega_1$. Choose a tree-like ladder system $\bar{\eta}=\left<
\eta_{\delta} : \delta \in S \right>$ such that
$\eta_{\delta}(\alpha)$ is a successor ordinal for all $\alpha <
\omega$ and $\delta \in S$. We enumerate the sets $P_{\alpha}$ by
$\omega$ without repetitions, e.g. $P_{\alpha}=\{ p_{\alpha, n} :
n < \omega \}$. Now let $F$ be the free group generated by the
elements $\{ x_{\nu} : \nu < \omega_1 \} \cup \{ y_{\delta,n} :
\delta \in S, n < \omega \}$. Let
$z_{\delta,-1}=y_{\delta,0}/p_{\delta,0}$ and for $n \geq 0$
\[ z_{\delta,n}=(y_{\delta,0} -
w_{\delta,n})/(\prod\limits_{i=0}^{n+1} p_{\delta,i}), \] where
$w_{\delta,n}=\sum\limits_{i=0}^n(\prod\limits_{j=0}^ip_{\delta,j})x_{\eta_{\delta}(i)}$.
Let $G$ be the subgroup of $F$ generated by the elements $\{
x_{\nu}:\nu < \omega_1 \} \cup \{z_{\delta,n} : \delta \in S, n <
\omega\}$. Then the only relations between the generators of $G$
are
\begin{equation}
\label{equ}  p_{\delta,n+1}z_{\delta,n}=z_{\delta,n-1} -
x_{\eta_{\delta}(n)}
\end{equation}
for $\delta \in S$ and $n \geq 0$. Now, for $\nu < \omega_1$ let
$G_{\nu}$ be the pure closure in $G$ of $G \cap \left(\{x_{\mu} :
\mu < \nu \} \cup \{ z_{\delta,n} : \delta \in S\cap \nu, n <
\omega \}\right)$. Then the sequence $\left< G_{\nu} : \nu <
\omega_1 \right>$ forms an $\omega_1$-filtration of $G$.
Moreover, for $\nu \in S$ we have \begin{equation} G_{\nu
+1}/G_{\nu} \cong F_{\nu} \oplus H_{\nu} ,\end{equation} where
$F_{\nu}$ is the free group on the generator $x_{\nu} + G_{\nu}$
and $H_{\nu} \cong \left< 1/p_{\nu,n} : n < \omega
\right>=:R_{P_{\nu}} \subseteq \Q$. Finally, $G$ is
$\aleph_1$-free by Pontryagin's criterion. Indeed, if $J_0$ is a
finite subset of $G$, then the pure closure of $J_0$ is contained
in the pure closure of a finite subset $J_1$ of $\{ x_{\nu} : \nu
< \omega_1 \} \cup \{ y_{\delta,n} : \delta \in S, n < \omega
\}$. By enlarging $J_1$ we may assume that there exists $m$ such
that for all $y_{\delta,o} \in J_1$, $x_{\eta_{\delta}(n)} \in
J_1$ if and only if $n \leq m$. Then the equations (\ref{equ})
show that the pure closure of $J_1$ is
free (compare \cite[Example VIII 1.1]{EM}).\\
Finally, let $P$ be an infinite subset of $\bar{\Pi}$. Since $G$ is
$\aleph_1$-free there exists no finite rank pure subgroup $H$ of
$G$ such that $T_P \not\in \tc(H)$, hence Proposition
\ref{characterization} shows that $T_P \not\in \tc(G)$ if and only
if the set \begin{equation} N=\{ \delta < \omega_1 : G/G_{\delta}
\textit{ contains a finite rank pure subgroup $L_{\delta}
\subseteq G/G_{\delta}$}\end{equation} \[\textit{ such that }
T_P \not\in \tc(L_{\delta}) \}
\] is stationary in $\omega_1$. Since for $\delta \in S$ we have
$G_{\delta+1}/G_{\delta}\cong R_{P_{\delta}}$ it is now easy to
see that $N$ is stationary if and only if $U=\{ \delta < \omega_1
: |P \cap P_{\delta}|=\aleph_0 \}$ is stationary. Note that $S$
is a cub in $\omega_1$.
\end{proof}

If $G$ is a torsion-free group, then it is not hard to see that the set
\begin{equation}
\{ P \subseteq \bar{\Pi} : T_P \in \tc(G) \}
\end{equation}
forms an ideal on $\P(\bar{\Pi})$ containing all finite subsets
of $\bar{\Pi}$. In fact, the next theorem shows that every such
ideal may appear. To avoid additional notation let us allow an
ideal in $\P(\bar{\Pi})$ to contain $\bar{\Pi}$ itself.

\begin{theorem}[$CH$]
\label{main} Let $I \subseteq \P(\bar{\Pi})$ be an ideal
containing all finite subsets of $\bar{\Pi}$. Then there exists
an $\aleph_1$-free group $G$ of cardinality $\aleph_1$ such that
for every $P \subseteq \bar{\Pi}$, $T_P \in \tc(G)$ if and only if $P \in I$.
\end{theorem}

\begin{proof}
Let $I$ be given. If $\bar{\Pi} \in I$, then we choose $G$ to be
free of cardinality $\aleph_1$ and we are done. Therefore, assume
that $\bar{\Pi}\not\in I$. Choose a continuous increasing
sequence of boolean subalgebras $\left< B_{\alpha} \subseteq
\P(\bar{\Pi}): \alpha < \omega_1 \right>$ such that each
$B_{\alpha}$ is countable and contains all finite subsets of
$\bar{\Pi}$. Note that this is possible since we are assuming
$CH$. Let $\alpha < \omega_1$ and put
\begin{equation}
I \cap B_{\alpha} = \{ I^-_{\alpha,i} : i < \omega \}
\end{equation}
and
\begin{equation}
B_{\alpha} \backslash I =\{ I_{\alpha,i}^+ : i < \omega \},
\end{equation}
where we assume that each $I_{\alpha,i}^+$ is repeated infinitely
many times. Choose for $\alpha < \omega_1$ and $i < \omega$
\begin{equation}
p_{\alpha,i} \in I_{\alpha,i}^+ \backslash \left( \bigcup \{
I_{\alpha,j}^- : j <i \} \cup \{ i_{\alpha,j} : j < i\}\right).
\end{equation}
Note, that this is possible since $I_{\alpha,i}^+$ is infinite and
$\bigcup \{ I_{\alpha,j}^- : j <i \} \in I$. Let
\begin{equation}
P_{\alpha} = \{ p_{\alpha,i} : i < \omega \}
\end{equation}
and let $G$ be the group from Proposition \ref{existence} for
$\left< P_{\alpha} : \alpha < \omega_1 \right>$. Then $G$ is
$\aleph_1$-free and of cardinality $\aleph_1$ and by Proposition
\ref{existence} it suffices to prove that for a subset $P
\subseteq \bar{\Pi}$ we have $P \in I$ if and only if there
exists $\gamma < \omega_1$ such that for all $\delta > \gamma$,
$P \cap P_{\alpha}$ is finite. Thus let $P \subseteq \bar{\Pi}$
and assume that $P \in I$. Then there exists $\gamma < \omega_1$
such that for all $\delta > \gamma$, $P \in I \cap B_{\delta}$.
Fix $\delta > \gamma$, then $P=I^-_{\delta,i}$ for some $i <
\omega$. Hence, for all $j > i$ we obtain $p_{\delta,j} \not\in
I_{\delta,i}^-$. Thus, $P \cap P_{\delta} \subseteq
\{p_{\delta,j} : j \leq i \}$ which is finite. Conversely, assume
that $P \not\in I$. Then there exists $\gamma < \omega_1$ such
that for all $\delta> \gamma$, $P \in B_{\delta}\backslash I$.
Fix $\delta > \gamma$, then $P=I_{\delta,j}^+$ for infinitely
many $j < \omega$ by the choice of the $I_{\delta,i}^+$'s. But, if
$P=I_{\delta,j}^+$, then $p_{\delta,j} \in P_{\delta} \cap
(P\backslash \{p_{\delta,i} : i <j \})$ and hence $P \cap
P_{\delta}$ is infinite. This finishes the proof.
\end{proof}

 We are now able to characterize the class
$\tc(G)$ for torsion-free groups of cardinality $\aleph_1$ assuming $CH$.

\section{The characterization}

In \cite[Theorem 3.6]{SW} a characterization of all classes $\C$
of torsion groups was given which could satisfy $\C=\tc(C)$ for
some completely decomposable group $C$. We shall show next that we
can drop condition \cite[Theorem 3.6 (v)]{SW} if we assume $CH$
and replace completely decomposable by $\aleph_1$-free of
cardinality $\aleph_1$. Recall, that condition \cite[Theorem 3.6
(v)]{SW} says the following
\begin{equation}
\label{*} \textit{ If $P$ is an infinite set of primes such that
$T_P \not\in \C$,}
\end{equation}
\[ \textit{ then there exists an infinite subset $P^{\prime}$ of $P$ such that} \]
\[ \textit{ for all infinite subsets $X$ of $P^{\prime}$, $T_X \not\in \C$.}
\]

\begin{theorem}[$CH$]
\label{cor} Let $\C$ be a class of torsion groups. Then
$\C=\tc(G)$ for some ($\aleph_1$-free) torsion-free group $G$ of
cardinality less than or equal to $\aleph_1$ if and only if the
following conditions are satisfied:
\begin{enumerate}
\item $\C$ is closed under epimorphic images;
\item $\C$ contains all torsion cotorsion groups;
\item If $p$ is a natural prime, then $\bigoplus\limits_{n < \omega}\Z(p^n) \in \C$ if and only if
$\C$ contains all $p$-groups;
\item If $P$ is an infinite subset of $\Pi$, then $\bigoplus\limits_{p \in P}\Z(p) \in \C$ if and only if
$\bigoplus\limits_{p \in P}H_p \in \C$ for all $p$-groups $H_p
\in \C$ ($p \in P$).
\end{enumerate}
\end{theorem}

\begin{proof}
Let us first show that (i) to (iv) hold for $\tc(G)$ for any
torsion-free group $G$ of cardinality less than or equal to
$\aleph_1$. Clearly, (i) and (ii) are true. Moreover, if $G$ is
countable, then \cite[Corollary 3.7]{SW} shows that (iii) and (iv)
hold for $G$. Thus assume that $G$ is of cardinality $\aleph_1$
and let $\left< G_{\alpha} : \alpha < \omega_1 \right>$ be an
$\omega_1$-filtration of $G$. Let $p$ be a prime and assume that
$\bigoplus\limits_{n < \omega}\Z(p^n) \in \tc(G)$. Moreover,
assume that $T$ is a $p$-group and $T \not\in \tc(G)$. By
Proposition \ref{characterization} there exists either a finite
rank pure subgroup $H$ of $G$ such that $T \not\in \tc(H)$ or the
set $Q=\{ \delta < \omega_1 : G/G_{\delta} \textit{ contains a
finite rank pure subgroup $L_{\delta} \subseteq G/G_{\delta}$ such
that } T \not\in \tc(L_{\delta}) \}$ is stationary in $\omega_1$.
If $H$ exists, then \cite[Theorem 3.6]{SW} shows that
$\bigoplus\limits_{n < \omega}\Z(p^n) \not\in \tc(H)$
contradicting the fact that $\bigoplus\limits_{n < \omega}\Z(p^n)
\in \tc(G) \subseteq \tc(H)$. Thus assume that $Q$ is stationary
in $\omega_1$. Again by \cite[Theorem 3.6]{SW} it follows that
for $\delta \in Q$ also $\bigoplus\limits_{n < \omega}\Z(p^n)
\not\in \tc(L_{\delta})$ since all $G_{\alpha}$'s are countable.
Thus \begin{equation} Q=\{ \delta < \omega_1 : G/G_{\delta}
\textit{ contains a finite rank pure subgroup $L_{\delta}
\subseteq G/G_{\delta}$}\end{equation} \[ \textit{ such that }
\bigoplus\limits_{n < \omega}\Z(p^n)\not\in \tc(L_{\delta}) \}\]
and Proposition \ref{characterization} shows that
$\bigoplus\limits_{n < \omega}\Z(p^n) \not\in \tc(G)$ - a
contradiction. Thus (iii) holds since the converse implication is
trivial.\\
It is straightforward to see that also (iv) holds for $\tc(G)$ using similar arguments as above.\\
Finally, assume that $\C$ satisfies (i) to (iv). We identify
$\omega$ with $\bar{\Pi}$ by a bijection $i: \omega \rightarrow
\bar{\Pi}$. Let $I=\{ X \subseteq \omega : \bigoplus\limits_{p \in
i(X)}\Z(p) \in \C \}$. Then it is easy to see that $I$ is an
ideal on $\omega$ containing all finite subsets of $\omega$. Thus
by Theorem \ref{main} there exists an $\aleph_1$-free group $G$ of
cardinality $\aleph_1$ such that for every subset $P \subseteq
\bar{\Pi}$, $\bigoplus\limits_{p \in P}\Z(p) \in \tc(G)$ if and
only if $i^{-1}(P) \in I$. Since we have already shown that
$\tc(G)$ satisfies (i) to (iv) it is now obvious that $\C=\tc(G)$.
\end{proof}

Since it was shown in \cite[Theorem 2.7]{St2} and \cite[Corollary
3.9]{SW} that in G\"odel's universe for every torsion-free group $G$
Theorem \ref{cor} (i) to (iv) are satisfied for $\C=\tc(G)$ we immediately get
the following result.

\begin{corollary}[$V=L$]
For every torsion-free group $G$ there exists an $\aleph_1$-free
group $H$ of cardinality $\aleph_1$ such that $\tc(G)=\tc(H)$.
\end{corollary}

Moreover, we obtain the existence of some ugly torsion-free
groups showing that the $\tc$-Conjecture from
\cite[$\tc$-Conjecture 2.12]{St2} does not hold. In \cite{St2} it
was conjectured that in $V=L$ for every torsion-free group $G$
there exists a completely decomposable group $C$ such that
$\tc(G)=\tc(C)$, hence condition (\ref{*}) would be satisfied for
all torsion-free groups $G$. This is not the case.

\begin{corollary}[CH]
For every infinite set of primes $P$ there exists an $\aleph_1$-free torsion-free
group $G$ of cardinality $\aleph_1$ satisfying
$T_P \not\in \tc(G)$ such that for every infinite subset $Q
\subseteq P$ there exists an infinite subset $Q_1 \subseteq Q$
such that $T_{Q_1} \in \tc(G)$. Thus $\tc(G) \not= \tc(C)$ for
every completely decomposable group $C$.
\end{corollary}

\begin{proof}
Let $P$ be the given infinite set of primes. It was shown by Eda
in \cite[Proof of Theorem 5]{E} that there exists a strictly
decreasing chain of subsets $P_{\alpha} \subseteq P$ ($\alpha <
\omega_1$) such that
\begin{enumerate}
\item $P_{\alpha}$ is infinite;
\item $\alpha < \beta$ implies $P_{\beta}$ is almost contained in $P_{\alpha}$;
\item $\alpha < \beta$ implies $|P_{\alpha}\backslash P_{\beta}|$ is infinite;
\item $\bigcap\limits_{\alpha < \omega_1}P_{\alpha}$ is finite.
\end{enumerate}
Let $U$ be the ultrafilter generated by $\bar{P}=\{ P_{\alpha} :
\alpha < \omega_1 \}$ and let $G$ be the group from Proposition
\ref{existence} for $\bar{P}$. If $Q$ is an infinite subset of
$P$, then divide $Q$ into two disjoint infinite subsets, e.g.
$Q=Q_1 \cup Q_2$. Since $U$ is an ultrafilter it follows that
without loss of generality $Q_1 \not\in U$. Hence, there exists
$\alpha < \omega_1$ such that $|P_{\alpha} \cap Q_1|$ is finite.
Thus, for every $\alpha \leq \beta$ we obtain $| P_{\beta} \cap
Q_1|$ is finite. Therefore, the set $\{ \delta < \omega_1 : |P
\cap P_{\delta}|=\aleph_0 \}$ is not stationary in $\omega_1$ and
Proposition \ref{existence} implies that $T_{Q_1}\in \tc(G)$.\\
Finally, $\tc(G)\not=\tc(C)$ for any completely decomposable group
$C$ since $\tc(G)$ violates \cite[Theorem 3.2 (v)]{SW} which is
our condition (\ref{*}).
\end{proof}

\goodbreak

\end{document}